\documentclass[12pt]{article}
\usepackage{sigsam, amsmath}

\issue{TBA}
\articlehead{TBA}
\titlehead{Algorithms for quadratic forms over global function fields of odd characteristic}
\authorhead{Mawunyo Kofi Darkey-Mensah}
\begin{document}

\title{Algorithms for quadratic forms over global function fields of odd characteristic}

\author{Mawunyo Kofi Darkey-Mensah \\
Institute of Mathematics \\
University of Silesia \\
Katowice, Poland, 40-007 \\
\url{mdarkeymensah@gmail.com}}

\date{}

\maketitle

\begin{abstract}
This paper presents an adaptation of recently developed algorithms for quadratic forms over number fields in \cite{koprowski2018computing} to global function fields of odd characteristics. First, we present algorithm for checking if a given non-degenerate quadratic form is isotropic or hyperbolic. Next we devise a method for computing the dimension of the anisotropic part of a quadratic form. Finally we present algorithms computing two field invariants: the level and the Pythagoras number.
\end{abstract}

\section{Introduction}
The theory of quadratic forms over fields is a well developed area of mathematics. However little work has been done so far for the computational aspects of the theory. Majority of previously published research focused mainly on forms over rationals. In a much recent work by P. Koprowski and A. Czogała in \cite{koprowski2018computing}, the authors developed a number of algorithms for quadratic forms over number fields. Namely, algorithms for checking the isotropy and hyperbolicity of forms, determining the dimension of an anisotropic form, and computing two field invariants, the level and Pythagoras number of a number field. The goal in this paper is to extend thier results to forms over global function fields of characteristic different from 2.

In this paper, the field $K$ is always a global function field of characteristic $\neq 2$ (thus a finite extension of the field of rational functions in one variable over a finite field), and $\OO_K$ will denote the integral closure of $\FFq[x]$ in $K$. The construction of majority of the algorithms in this work reflects the idea of solving the problem locally, in a similar order and fashion as present in \cite{koprowski2018computing}. The local solutions are then used as sub-procedures to find the global solutions.

The paper is organized as follows: in Sections \ref{ch3_sec_isotropy} we present algorithms to check if a given quadratic form over a global function field is isotropic or hyperbolic. Next in Section \ref{ch3_sec_witt_index}, Algorithm \ref{anisotropic_dimension_global} computes the anisotropic dimension of a form in a global function field. This algorithm is important for computing the Witt index of a quadratic form. In chapter \ref{lengthofsos} we present an algorithm that finds a
minimal number of squares needed to represent a given sum of squares in a global function field. This is known as the length of a sum of squares. Finally in Sections \ref{ch3_sec_level}, we present algorithms for computing two important field invariants, namely the level $s(K)$ and the Pythagoras number $P(K)$ of a global function field.




\section{Isotropy and Hyperbolicity of a quadratic form}
\label{ch3_sec_isotropy}
Recall that a quadratic form $q$ is said to be isotropic over a field $K$ if there exists a non-zero vector $v$ with coefficients in $K$ such that $q(v)=0$. Observe that over a finite field, an isotropy test is trivial. It suffices to look at it's discriminant. A binary form over a residue field is isotropic if and only if it's determinant is a minus square (see. \cite[Theorem I.3.2]{lam2005introduction}). Finally, forms of $\dim \geq 3$ over finite fields are always isotropic by means of \cite[Example XI.6.2]{lam2005introduction}. For local fields, the problem is similarly easy. A unary form is never isotropic and quintic or higher dimensional forms are always isotropic by means of \cite[Theorem VI.2.12]{lam2005introduction}. Suppose that $q \cong q_0 \perp \Ideal{\pi}q_1$ where $q_0 = \Ideal{a_1,\dots,a_r}$ and $q_1 = \Ideal{a_{r+1},\dots,a_n}$ with all $a_i$ being $\gp$-adic units. Then by \cite[Theorem VI.1.9]{lam2005introduction} $q$ is isotropic if either $q_0$ or $q_1$ is isotropic over the residue field.  Thus, over local fields isotropy is easy to check. We use this fact to check isotropy over a global function field.\\

\begin{algorithm}[H]\caption{Isotropy in a global function field $K$.}
\label{isotropy_global}
\KwIn{A non-degenerate diagonal quadratic form $q = \Ideal{a_1,\dots,a_n}$ over $K$ with $a_i \in \OO_K$}
\KwOut{True if $q$ is isotropic and False otherwise.}
\If{$\dim q \leq 1$}{\Return False}
\If{$\dim q=2$}{
    \eIf{$\disc q$ is a square in $K$}{
        \Return True}{
        \Return False}}
Let $\gP = \{\gp_1,\dots,\gp_n\}$ be the list of all places of $K$ dividing any of the coefficients $a_i$ of $q$\;
\For{$\gp \in \gP$}{
    \If{$q$ is not isotropic over $K_\gp$}{
        \Return False}}
\Return True
\end{algorithm}

\begin{poc}
A unary form is always anisotropic, and it is well known (see. e.g. \cite[Theorem I.3.2]{lam2005introduction}) that a binary form $q$ is isotropic if $\disc q$ is a square, and anisotropic otherwise. For forms of higher dimension we use the local-global principle \cite[Principle VI.3.1]{lam2005introduction}. The form is isotropic over $K$ if and only if it is isotropic over all the completions of $K$. Now $q$, having dimension at least three, is trivially isotropic at all places that do not divide any of the coefficients. The remaining places to check are all the places dividing any of the coefficients in $q$. 
\end{poc}

Next, a quadratic form is hyperbolic if it is an orthogonal sum of binary forms isometric to $\Ideal{1,-1}$. A hyperbolic form is universal in the sense that it represents all the zero elements of the field (see e.g. \cite[I.3]{lam2005introduction}). We first observe from \cite[Theorem I.3.2]{lam2005introduction} that if $q$ is hyperbolic, then it is even dimensional. Furthermore, if the discriminant $\disc q$ is a square in $K_\gp$ and the Hasse invariant $h_\gp(q)$ equals $(-1,-1)^{m(m-1)/2}_\gp$, then $q$ is isometric to the hyperbolic space $m\Ideal{1,-1}$ by \cite[Proposition V.3.25]{lam2005introduction} and hence hyperbolic, otherwise $q$ is not hyperbolic over the local field $K_\gp$. Hence like with isotropy, it is easy to check if a form is locally hyperbolic. We again use it to check if a form is globally hyperbolic. \\

\begin{algorithm}[H]\caption{Hyperbolicity in a global function field $K$.}
\label{hyperbolicity}
\KwIn{A non-degenerate diagonal quadratic form $q=\Ideal{a_1,\dots,a_n}$ over a global function field $K$.}
\KwOut{True if $q$ is hyperbolic and False otherwise.}
\If{$\dim q$ is odd}{\Return False}
\If{$\disc q$ is not a square in $K$}{\Return False}
Let $\gP =\{\gp_1,\dots,\gp_m\}$ be the list of places dividing any of the coefficients $a_1,\dots,a_n$ of $q$ in $\OO_K$\;
\For{$\gp \in \gP$}{
\tcp{Check if $q\otimes K_\gp$ is hyperbolic.}
    \If{$q$ is not hyperbolic over $K_\gp$}{
    \Return False} }
\Return True
\end{algorithm}

\begin{poc}
It is well known (see e.g. \cite[Theorem I.3.2]{lam2005introduction}) that the discriminant of a hyperbolic form is a square and it's dimension has to be even. Moreover, by the Weak Hasse-Minkowski Principle \cite[Corollary VI.3.3]{lam2005introduction}, a quadratic form is hyperbolic over a global field if and only if it is hyperbolic over every completion of the field. The places of importance in this situation are those dividing any of the coefficients of the form $q$.
\end{poc}

\section{Witt index of a quadratic form}
\label{ch3_sec_witt_index}
Recall (see e.g. \cite[Theorem I.4.1]{lam2005introduction}) that any non-degenerate quadratic form $q$ can be uniquely decomposed as $q = q_h \perp q_a$, where $q_h$ is hyperbolic (or zero), and $q_a$ is anisotropic. The Witt index of $q$, denoted $\ind(q)$, is the number of hyperbolic planes constituting $q_h$, i.e. half of the dimension of $q_h$. It can be computed with the formula $\ind(q) = \frac{1}{2} \cdot (\dim q - \dim q_a)$. In the following algorithm we make use of \cite[Algorithm 8]{koprowski2018computing} to compute the anisotropic dimension of the localization of $q$ at some place $\gp$ of $K$. Although \cite[Algorithm 8]{koprowski2018computing} was designed for completions of number fields, it can be used without any modifications in our case as well.

\begin{algorithm}[H]\caption{Anisotropic dimension in a global function field $K$} \label{anisotropic_dimension_global}
\KwIn{A non-degenerate diagonal quadratic form $q=\Ideal{a_1,\dots,a_n}$ over $K$.}
\KwOut{Dimension of the anisotropic part of $q$.}
$Dimensions \gets [ \ ]$\;
\For{$a \in q$}{
Let $\gP = \{\gp_1,\dots,\gp_m\}$ be a list of places dividing $a$ in $\OO_K$\;
    \For{$\gp \in \gP$}{
    \tcp{Apply Algorithm \cite[Algorithm 8]{koprowski2018computing} to compute the dimension of the anisotropic part of $q \otimes K_\gp$}
    $d_\gp \gets \texttt{AnisotropicDimLocal($q,\gp$)}$\;
    Append $d_\gp$ to $Dimensions$\;
    }
}
\Return $\max Dimensions$
\end{algorithm}

\begin{poc}
The dimension of the anisotropic part of $q$ is clearly the maximum of the dimensions of the anisotropic parts of the localizations of $q$ at places of $K$ which divide any coefficient in $q$.
\end{poc}

\section{Length of an element}\label{lengthofsos}
The length of an element $a \in K$, denoted $\len(a)$, is the smallest natural number $n \in \NN$ such that $a$ can be written as a sum of $n$ squares in $K$. In order to determine the global length of an element $a$, we need to first compute the local lengths at all places of $K$ dividing $a$. Recall that if $\gp$ is a place of $K$, then the square class group of the local field $K_\gp$ has the form $\dot{K}_\gp/\dot{K}^2_\gp = \{ \dot{K}^2_\gp, u_\gp \dot{K}^2_\gp, \pi_\gp \dot{K}^2_\gp, u_\gp \pi_\gp \dot{K}^2_\gp\}$ where $\ord_\gp u_\gp \equiv 0 \modd{2}$ is a $\gp$-adic unit, and $\ord_\gp \pi_\gp \equiv 1 \modd{2}$ is a $\gp$-adic uniformizer. Now if the $\gp$-adic valuation $v_\gp(a)$ is even. Then either $a=1$ or $a=u_\gp$ (mod squares). If $a$ is a square, then the local length $\len_\gp(a)=1$. Otherwise, the Hilbert symbol $(-1,a)_\gp = (-1,u_\gp)_\gp = 1$, i.e. $1 \in D_\gp(\langle -1,a \rangle) \cong a \in D_\gp(\langle 1,1 \rangle)$, hence $\len_\gp(a)=2$. If $v_\gp(a)$ is odd. Then either $a=\pi_\gp$ or $u_\gp \pi_\gp$ (mod squares) and hence $(-1,a)_\gp = (-1,u_\gp \pi_\gp)_\gp = (-1,\pi_\gp)_\gp$ is either $1$ if $-1 \in \dot{K}_\gp^2$, or $-1$ if $-1 \notin \dot{K}_\gp^2$. If $-1 \in \dot{K}_\gp^2$ ($\cong$ $\langle 1,1 \rangle$ isotropic over $K(\gp)$), then $\len_\gp(a)=2$. If $-1 \notin \dot{K}_\gp^2$, then the form $\langle 1,1,1 \rangle$ is isotropic over $K_\gp$. Hence $a \in D_\gp(\langle 1,1,1 \rangle)$ and $\len_\gp(a)=3$. Similarly as in the previous sections, we make use of the length over a local field to compute the length over a global function field.\\

\begin{algorithm}[H] \caption{Length in a global field} \label{globallength}
\KwIn{A nonzero element $a$ of a global function field $K$}
\KwOut{Length of $a$ in $K$}
\eIf{$a$ is a square in $K$}{
    \Return $1$}
{
$ \gL \gets [ 2 ]$\;
Let $\gP=\{ \gq_1, \dots, \gq_n\}$ be the list of places dividing $a$ in $K$\;
\For{$\gq \in \gP$}{
    Compute the length $\len_\gq(a)$ of $a$ in the completion $K_\gq$\;
    \If{$\len_\gq(a) = 3$}{
        \Return $3$}
    Append $\len_\gq(a)$ to $\gL$\;
}
\Return max $\gL$\;
}
\end{algorithm}

\begin{remark}
Let us mention that Algorithm \ref{globallength} is a part of the joint paper \cite{darkey2021computing} of the author and Beata Rothkegel, which is presently under review elsewhere.
\end{remark}

\section{Level  and Pythagoras number}
\label{ch3_sec_level}
In this section, we present two important field invariants of a global function field, namely the level and Pythagoras number. Recall (see e.g. \cite[\S XI.2 \& \S XI.5]{lam2005introduction}) that the level of a field $K$ (denoted $s(K)$) is the \emph{length} of $-1$ in $K$, and a Pythagoras number (denoted $P(K)$) of a field $K$ is the smallest positive integer $n \in \NN$ such that every sum of squares in $K$ is a sum of $n$ squares. Below Algorithms \ref{level} and \ref{Pythagoras_number} computes the level and Pythagoras number, respectively. \\

\begin{algorithm}[H]\caption{Level}\label{level}
\KwIn{A global function field $K$ with full field of constants $\FFq$ of order $q$.}
\KwOut{The level $s(K)$ of $K$.}
\If{$q \equiv 1 \modd{4}$}{\Return 1}
\If{$q \equiv 3 \modd{4}$}{\Return 2}
\end{algorithm}

\bigskip

\begin{algorithm}[H]\caption{Pythagoras number}\label{Pythagoras_number}
\KwIn{A global function field $K$ with full field of constants $\FFq$ of order $q$.}
\KwOut{The Pythagoras number $P(K)$ of $K$.}
\If{$q \equiv 1 \modd{4}$}{\Return 2}
\If{$q \equiv 3 \modd{4}$}{\Return 3}
\end{algorithm}

\bigskip

The algorithms computing the level and Pythagoras number in a global function field are very similar, hence we prove the correctness of Algorithms \ref{level} and \ref{Pythagoras_number} together below.

\begin{poc}
If $q \equiv 1 \modd{4}$, then $-1 \in \dot{\FFq}^2 \subset \dot{K}^2$, so $s(K)=1$. This implies that the form $\langle 1,1 \rangle$ is isotropic over $K(\gp)$ for all $\gp \in \Omega(K)$. We thus have $P(K)=2$. 

Conversely, if $q \equiv 3 \modd{4}$, then $-1 \notin \dot{K}^2$; but $-1 \in D_K(2)$ by means of \cite[Proposition II.3.4]{lam2005introduction}, so $s(K)=2$. The form $\langle 1,1,1 \rangle$ is isotropic over $K_\gp$ for all $\gp \in \Omega(K)$ and consequently $P(K)=3$.
\end{poc}

The presented algorithms can be implemented in existing computer algebra systems. Indeed, one can find a recent implementation in CQF – a free, open-source Magma \cite{magma} package for doing computations in quadratic forms theory (see \cite{koprowski2020cqf}).


\bibliography{ref}

\end{document}